\documentclass{gtart_h}

\def\ifplaintex{\expandafter\ifx\csname documentclass\endcsname\relax}

\def\gtp{{\mathsurround=0pt\it $\cal G\mskip-2mu$eometry \&\ 
$\cal T\!\!$opology $\cal P\!$ublications}}  

\def\recd{{\small Received:\qua\receiveddate\ifx\reviseddate\relax
\else\qquad Revised:\qua\reviseddate\fi\par}} 

\def\lognumber#1{\def\thelognumber{#1}}
\def\volumenumber#1{\def\thevolumenumber{#1}}
\def\volumeyear#1{\def\thevolumeyear{#1}}
\def\papernumber#1{\def\thepapernumber{#1}}
\def\pagenumbers#1#2{\def\startpage{#1}\def\finishpage{#2}}
\def\published#1{\def\publishdate{#1}}

\def\received#1{\def\receiveddate{#1}}

\def\accepted#1{\def\accepteddate{#1}}

\long\def\asciiabstract#1{\long\def\theasciiabstract{#1}}

\let\\\par\let\thelognumber\relax\let\thevolumenumber\relax
\let\thepapernumber\relax\let\thevolumeyear\relax\let\startpage\relax
\let\finishpage\relax\let\publishdate\relax\let\receiveddate\relax
\let\reviseddate\relax\let\accepteddate\relax\let\theasciititle\relax
\let\theasciiauthors\relax
\let\theasciiabstract\relax

\let\theasciiemail\relax

\ifplaintex
\font\logobig=cmssbx10 scaled 3836
\font\logomed=cmssbx10 scaled 2557
\else
\font\logobig=cmssbx10 scaled 4200
\font\logomed=cmssbx10 scaled 2800
\fi

\long\def\makeagttitle{   
\count0=\startpage
\agt\hfill      
\hbox to 45truept{\vbox to 0pt{\vglue -13truept{\logomed A\kern -.37em{\logobig 
T}\kern -.38em G}\vss}\hss}
\break
{\small Volume \thevolumenumber\ (\thevolumeyear)
\startpage--\finishpage\nl
Published: \publishdate}

\vglue .25truein

{\parskip=0pt\leftskip 0pt plus
1fil\def\\{\par\smallskip}{\Large\bf\thetitle}\par\medskip} \vglue
0.05truein

{\parskip=0pt\leftskip 0pt plus 1fil\def\\{\par}{\sc\theauthors}
\par\medskip}%
 
\vglue 0.03truein 

{\small\leftskip 25truept\rightskip 25truept{\bf Abstract}\stdspace\theabstract

{\bf AMS Classification}\stdspace\theprimaryclass
\ifx\thesecondaryclass\relax\else; \thesecondaryclass\fi\par
{\bf Keywords}\stdspace \thekeywords\par}\vglue 7truept

}   

\ifplaintex
\font\phead=cmsl9 scaled 950
\font\pnum=cmbx10 scaled 913
\font\pfoot=cmsl9 scaled 950
\headline{\vbox to 0pt{\vskip -4.5mm\line{\small\phead\ifnum
\count0=\startpage ISSN 1472-2739 (on-line) 1472-2747 (printed)
\hfill {\pnum\folio}\else\ifodd\count0\def\\{ }%
\ifx\theshorttitle\relax\thetitle\else\theshorttitle\fi\hfill{\pnum\folio}
\else\def\\{ and }{\pnum\folio}\hfill\ifx\theshortauthors\relax\theauthors
\else\theshortauthors\fi\fi\fi}\vss}}
\footline{\vbox to 0pt{\vglue 0mm\line{\small\pfoot\ifnum\count0=\startpage
\copyright\ \gtp\hfill\else
\agt, Volume \thevolumenumber\ (\thevolumeyear)\hfill\fi}\vss}}
\else
\font=cmsl9 scaled 1050
\font\lnum=cmbx10 
\font=cmsl9 scaled 1050
\makeatletter
\def\@oddhead{{\small\lhead\ifnum\count0=\startpage ISSN 1472-2739 
(on-line) 1472-2747 (printed)\hfill {\lnum\number\count0}\else\ifodd\count0
\def\\{ }\ifx\theshorttitle\relax \thetitle \else\theshorttitle\fi\hfill
{\lnum\number\count0}\else\def\\{ and }{\lnum\number\count0}
\hfill\ifx\theshortauthors\relax 
\theauthors\else\theshortauthors\fi\fi\fi}}\def\@evenhead{\@oddhead}
\def\@oddfoot{\small\lfoot\ifnum\count0=\startpage\copyright\ \gtp\hfill\else
\agt, Volume \thevolumenumber\ (\thevolumeyear)\hfill\fi}
\def\@evenfoot{\@oddfoot}
\makeatother
\fi
\let\maketitlepage\makeagttitle
\let\makeshorttitle\maketitlepage
\let\maketitle\maketitlepage

\newwrite\gtoutfile
\long\gdef\makeheadfile{  
{\def\\{, }\def\s{ }
\immediate\openout\gtoutfile head.xxx
\immediate\write\gtoutfile{Proxy-for: \ifx\theasciiauthors\relax
\theauthors\else\theasciiauthors\fi\s<\ifx\theasciiemail\relax\theemail\else\theasciiemail\fi>}
\immediate\write\gtoutfile{\noexpand\\}
\immediate\write\gtoutfile{Authors: \ifx\theasciiauthors\relax
\theauthors\else\theasciiauthors\fi}
{\def\\{ }\immediate\write\gtoutfile{Title: \ifx\theasciititle\relax
\thetitle\else\theasciititle\fi}}
\immediate\write\gtoutfile{Subj-class: GT or SG, GR etc}
\immediate\write\gtoutfile{MSC-class: \theprimaryclass\ifx\thesecondaryclass\relax\else, \thesecondaryclass\fi}
\immediate\write\gtoutfile{Journal-ref: Algebr. Geom. Topol. \thevolumenumber\s
(\thevolumeyear) \startpage-\finishpage}
\immediate\write\gtoutfile{Comments: Published by Algebraic and
Geometric Topology at}
\immediate\write\gtoutfile{\s\s\s  http://www.maths.warwick.ac.uk/agt/AGTVol\thevolumenumber/agt-\thevolumenumber-\thepapernumber.abs.html}
\immediate\write\gtoutfile{\noexpand\\}
\immediate\write\gtoutfile{}
\ifx\theasciiabstract\relax
\immediate\write\gtoutfile{\theabstract}\else
\immediate\write\gtoutfile{\theasciiabstract}\fi
\immediate\write\gtoutfile{}
\immediate\write\gtoutfile{\noexpand\\}
\immediate\write\gtoutfile{}
\immediate\closeout\gtoutfile}}  

\def\maketitlepage{\makeagttitle\makeheadfile}
\let\makeshorttitle\maketitlepage
\let\maketitle\maketitlepage

\lognumber{41}
\volumenumber{4}
\volumeyear{2004}
\papernumber{41}
\pagenumbers{943}{960}
\received{17 February 2004} 
\accepted{21 September 2004}
\published{22 October 2004}

\usepackage{amssymb, amsmath}

\let\proc\begin
\def\items{\begin{enumerate}}
\def\enditems{\end{enumerate}}
\let\endprf\endproof

\newtheorem{Theorem}{Theorem}[section]

\theoremstyle{definition}

\swapnumbers
\newtheorem{Subsec}[Theorem]{}

\def\subsection#1{\goodbreak\begin{Subsec}\label{\theTheorem}
\addcontentsline{toc}{subsection}{\theTheorem #1}{$\phantom{.}$\bf #1}
\end{Subsec}\nobreak}

\def\calj{\cal J}

\def\goths{\mathfrak S\/}
\def\st{\mathop{\, \rm st}} 
\def\ost{\mathop{\, \rm \overline{st}\,}}

\def\hatLie{{\widehat{\,\Lie\,}}\>\!}
\def\Lie{\mathop{\rm Lie}}
\def\Hom{\mathop{\rm Hom}}
\def\Res{\mathop{\rm Res}}
\def\Ind{\mathop{\rm Ind}}
\def\lra{\longrightarrow}
\def\dee{\partial}
\def\hpi{\varphi}

\def\eps{\varepsilon}
\def\La{\Lambda}
\def\Tha{\Theta}
\def\tha{\theta}

\def\Om{\Omega}

\def\ga{\gamma}
\def\al{\alpha}
\def\la{\lambda}
\def\tha{\theta}
\def\sgm{\sigma}
\def\Sgm{\Sigma}

\def\bbbr{\mathbb R}
\def\bbbz{\mathbb Z}
\def\lan{\langle}
\def\ran{\rangle}

\def\til{\tilde}

\def\calc{\cal C}

\title{Partition complexes, duality and integral\\tree representations}

\authors{Alan Robinson}                  
\address{Mathematics Institute, University of Warwick, Coventry CV4 7AL, UK}

\email{car@maths.warwick.ac.uk}

\begin{abstract}
We show that the poset of non-trivial partitions 
of $\{1,2,\ldots,n\}$ has a fundamental 
homology class with coefficients in a Lie superalgebra.  
Homological duality then rapidly yields a range of known 
results concerning the integral representations of the 
symmetric groups $\Sgm_n$ and $\Sgm_{n+1}$ on the homology and 
cohomology of this partially-ordered set. 
\end{abstract}

\asciiabstract{%
We show that the poset of non-trivial partitions of 1,2,...,n has a
fundamental homology class with coefficients in a Lie superalgebra.
Homological duality then rapidly yields a range of known results
concerning the integral representations of the symmetric groups S_n
and S_{n+1} on the homology and cohomology of this partially-ordered
set.}


\begin{document}
\maketitle

\section{Introduction}\label{1}

This paper reveals the geometry underlying certain integral representations of the 
symmetric groups.  It claims few original results, but gives a unified geometrical 
and homological treatment and frequently strengthens the known theorems.  
Above all it aims to show that topological methods give efficient proofs of 
many of them.  

Geometric realization converts the lattice of non-trivial 
partitions of a set with $n$ elements into 
the space of fully-grown (that is, non-degenerate) 
trees with leaves labelled by the set $\{0,1,\dots,n\}$.  
This space has reduced homology in degree~$n-3$ only.  
Its integral cohomology is a representation of $\Sigma_n$ which is isomorphic to 
the Lie representation $\Lie_n$, 
up to a  twist by the sign character.  It follows from the 
$\Sigma_{n+1}$-symmetry of the tree space 
that $\Lie_n$ extends to an integral representation of~$\Sigma_{n+1}$. 

Many mathematicians have contributed to the identification of the cohomology of 
the partition lattice with 
the twisted Lie representation: excellent summaries of the history of the theory are 
given in \cite{Fre} and \cite{Wac}. 
We show that the isomorphism arises at the cocycle level 
from a duality structure: the partition lattice has a fundamental homology cycle 
with coefficients in a certain $\Sgm_n$-module~${\goths}_n$, and 
the isomorphism is induced by cap product.  The module ${\goths}_n$ is 
the multilinear part of a free Lie superalgebra on $n$ generators of odd degree: 
it is isomorphic to the twisted Lie representation, though not uniquely so. 

In the final section we give a geometrical derivation of a short exact sequence 
relating these representations, due to S. Whitehouse \cite{Whi}.  

The tree spaces are related to configuration spaces, and the
occurrence of the Lie representations in this topological situation
was observed by F.~R.~Cohen~\cite{Co1},\cite{Co2}.

\section{The space of fully-grown $n$-trees }\label{2}

From \cite{R-W} we recall topological properties of the space of $n$-trees, and 
prove 
a homeomorphism with the nerve of the lattice of partitions of the set 
$\{1,2,\dots,n\}$. 

A {\it tree\/} is a contractible 1-dimensional polyhedron $X$, and we shall 
always require $X$ to be compact.  There is a unique coarsest 
triangulation of $X$, in which no vertex lies on exactly two edges.  
A vertex meeting more than two edges is called a {\it node\/};  
a vertex incident upon only one edge is a {\it free vertex\/}.  

An edge which connects two nodes is an {\it internal \/} edge.  An edge with a free 
vertex is called a {\it leaf\/}.

Following \cite{Bor},\cite{B-V},  we introduce a moduli space of trees with a 
fixed number of leaves. 

\proc{Definition}\label{2.1}
An $n$-{\it tree\/} is a tree with the following extra structure.  
   \items
   \item Every internal edge $e$ is assigned a length $l_e$ where $0 < l_e \le 
1$. 
 All leaves are conventionally assigned a length of 1.  
   \item The free vertices are labelled by a bijective correspondence with the 
set 
$\{0,1,\dots,n\}$. 
   \enditems
\end{Definition}
If $n \ge 2$, then the labelling of free vertices is equivalent to a labelling of the 
leaves 
by a bijection with the same set $\{0,1,\dots,n\}$.  The leaf labelled 0 is also called 
the~{\it root\/}. 
Two $n$-trees are regarded as equivalent if there is an isometry between them 
which 
preserves the labelling.  A tree is {\it fully grown\/} if at least one internal edge has 
length 1. 

\proc{Example}\rm: 
{\bf the vertices $v_A$ of $\,T_n$}\label{2.2} 

Let $A$ be any subset of the set 
$[n] = \{0,1,\dots,n\}$ of labels, such that $A$ and its 
complement $[n] \setminus A$ each have more than one element.  

We define $v_A$ to be the fully grown $n$-tree which has just one internal edge 
$e_A$ 
(which necessarily has length 1), and which has the leaves labelled by elements 
of $A$ attached at one end of this, and the leaves labelled by the elements of 
$[n] \setminus A$ at the other:  

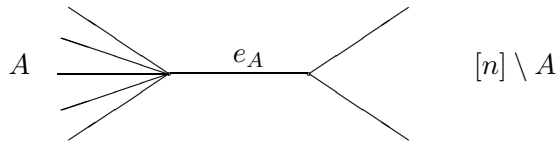
\setlength{\unitlength}{0.44pt}
\begin{picture}(450,170)(-375,-75)

\put(-59,0.4){\line(1,0){118}}

\put(-61,0.5){\line(-3,2){85}} 
\put(-61,0.3){\line(-3,1){92}} 
\put(-61,0){\line(-3,0){95}} 
\put(-61,-0.3){\line(-3,-1){92}}
\put(-61,-0.5){\line(-3,-2){85}}

\put(61,0.5){\line(3,2){85}}

\put(61,-0.5){\line(3,-2){85}}

\put(60,0){\circle{2}}
\put(-60,0){\circle{2}}

\put(-7,10){$e_A$}
\put(-200,0){$A$}
\put(200,0){$[n] \setminus A$}
\end{picture}
\end{Example}

Since trees are regarded as the same when they are related by a label-preserving 
homeomorphism, it follows that $v_A = v_{[n]\setminus A}$. 

In the next proposition $\Sigma_{n+1}$ denotes the set of all permutations of 
$\{0,1,\dots,n\}$, and 
$\Sigma_n$ the subgroup leaving 0 fixed. 

\proc{Proposition}\label{2.3} $\phantom{99}$
\items
\item The space $T_n$ of all fully-grown $n$-trees is a simplicial complex of 
dimension $n-3$, in 
which the vertices are the trees $v_A$ of {\/ \rm \ref{2.2}}. 
\item The space $\tilde{T}_n$ of all $n$-trees is the cone on $T_n$. 
\item The symmetric group $\Sigma_{n+1}$ acts on the right upon the pair 
$(\tilde{T}_n, T_n)$ by acting on the left upon the set of leaf labels.  This action is 
simplicial. 
\enditems\end{Proposition}

\proof  The proof is given in \cite{R-W}.  Suffice it to say here that $\tilde 
T_n$ 
has an evident cubical structure, with the internal edge-lengths as cubical 
coordinates. 
On the subspace $T_n$ the greatest edge-length is 1, so there is a simplicial 
structure 
with barycentric coordinates proportional to the edge-lengths: the vertices are the  
$v_A$ of \ref{2.2}, and each $k$-simplex is an amalgam of $k+1$ of the original cubes 
as 
illustrated for $k=2$ on the left in the following diagram. 

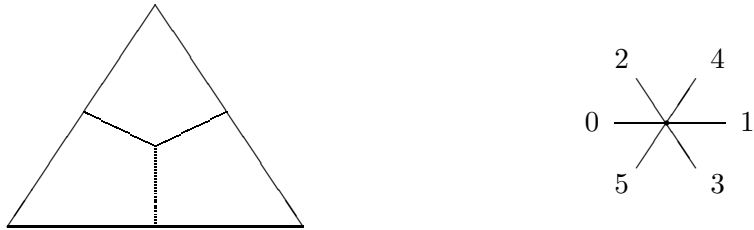
\setlength{\unitlength}{0.28pt}
\begin{picture}(450,500)(-325,-300)


\put(0,100){\line(2,-3){200}}
\put(0,100){\line(-2,-3){200}}
\put(-200,-200){\line(1,0){400}}


\multiput(0,-92)(0,-4){28}{\line(0,1){1}}
\multiput(0,-92)(3,1.45){33}{\line(0,1){1}}
\multiput(0,-92)(-3,1.45){33}{\line(0,1){1}}


\put(650,-120){\line(2,3){80}}
\put(650,0){\line(2,-3){80}}
\put(620,-60){\line(1,0){150}} 

\put(691,-60){\circle*{6}}


\put(580,-70){0}     
\put(790,-70){1}
\put(620,15){2}
\put(750,15){4}
\put(750,-155){3}
\put(620,-155){5}

\end{picture}

The space $\tilde T_n$ is evidently a cone with $T_n$ as its base. Its apex is a  
star-tree as illustrated on the right above, and the cone parameter maps  
 a tree to its maximal internal edge-length.  Therefore  $\tilde T_n$ can be 
 triangulated as the cone on the simplicial complex $T_n$.  As with  
$T_n$, the original cubical structure is a subdivision of this triangulation.  
\endprf

From now on, we use the triangulations just constructed, unless the contrary is 
mentioned. 

\proc{Lemma}\label{2.4} Every $(n-4)$-simplex of $\/ T_n$ is a face of exactly three 
top-dimensional simplices. \end{Lemma}

\proof  The top-dimensional simplices correspond to {\it binary trees \/} in 
which each node is the meet of three edges.  An $(n-4)$-simplex of $\/T_n$ 
arises from a tree with one exceptional node of order four: it faces three 
top-dimensional simplices corresponding to the three ways of resolving this node 
into two nodes of order three \cite{R-W}.
\endprf

Lemma \ref{2.4} shows that $T_n$ has three-to-one incidence on codimension-one 
simplexes.  As is well known, a triangulated closed manifold has two-to-one 
incidence, and consequently has   
a fundamental top-dimensional homology class (with integer coefficients, possibly 
twisted). We shall see in \S\ref{4} that $T_n$ too has a fundamental homology class for 
which the coefficients are a certain module~${\goths}_n$.

\proc{Theorem}\label{2.5} The space $T_n$ of fully-grown $n$-trees has the 
homotopy type of a wedge of $(n-1)!$ spheres of dimension $n-3$.  
The $\Sgm_n$-module $\tilde{H}_{n-3}(T_n)$  
restricts to the regular integral representation of the subgroup $\Sgm_{n-1}$. 
\end{Theorem}

\proof  It is shown in \cite{R-W} that $T_n$ is obtained from a contractible 
space by 
attaching $(n-1)!$ simplices of dimension $n-3$ along their boundaries.  These 
simplices correspond to the trees 

\setlength{\unitlength}{0.55pt}
\begin{picture}(450,170)(-325,-125)

\put(-200,0){\line(1,0){80}}
\put(-120,0){\line(1,0){80}} 
\put(-40,0){\line(1,0){80}}
\put(40,0){\line(1,0){30}}
\multiput(70,0)(6,0){10}{\line(1,0){1.5}} 
\put(130,0){\line(1,0){30}}
\put(160,0){\line(1,0){80}}

\put(-120,-0.5){\line(0,-1){75}}
\put(-40,-0.5){\line(0,-1){75}}
\put(41,-0.5){\line(0,-1){75}}
\put(160,-0.5){\line(0,-1){75}}

\put(-120,0){\circle{2}}
\put(-40,0){\circle{2}}
\put(40,0){\circle{2}}
\put(160,0){\circle{2}}


\put(-220,-8){$0$}
\put(-130,-98){$\sgm(1)$}
\put(-50,-98){$\sgm(2)$}
\put(30,-98){$\sgm(3)$}
\put(120,-98){$\sgm(n-1)$}
\put(250,-8){$n$}
\end{picture}

where $\sgm$ is any element of the permutation group $\Sgm_{n-1}$.  Thus 
$T_n$ has the homotopy 
type of a wedge of spheres, as claimed; and the homology classes of 
these spheres are 
regularly permuted by the subgroup $\Sgm_{n-1}$.  \endprf
 
\subsection{Homeomorphism with the nerve of the partition 
lattice $\La_n$}

We call a partition of the set $\{1,2,\dots,n\}$ 
{\it non-trivial\/} if it is neither indiscrete nor discrete: that is, the number of 
equivalence classes is greater than $1$ but less than $n$. The non-trivial 
partitions, 
ordered by refinement, form a lattice.  We use the notation $\La_n$ both for 
this lattice and for the simplicial complex which is its nerve.  

\proc{Proposition}\label{2.7} The nerve $\La_n$ of the partition lattice is 
$\Sgm_n$-equivariantly 
homeomorphic to the space $T_n$ of fully-grown $n$-trees. \end{Proposition}

\proof  We describe an explicit homeomorphism $\hpi\co T_n \to \La_n$.  A 
point $\al\in T_n$ is an $n$-tree in which each internal edge has a given length 
between $0$ and $1$. 
 
For each $1\le i \le n$  there is a unique arc $\ga_i$ in the tree $\al$ 
which starts at the root labelled $0$ and ends at the leaf labelled $i$. 
We parametrize all these arcs with unit speed, and we extend $\ga_i$ to a map 
$[0,\infty) \to \al$ by setting $\ga_i(t)$ to be constant at the vertex labelled $i$ for 
all $t\gg 0$. 

As time increases, these arcs diverge and determine finer and finer partitions of 
$\{1,2,\dots,n\}$. 
For every $t>0$ we define $\Pi_\al(t)$ to be the partition of $\{1,2,\dots,n\}$ given 
by 
the 
equivalence relation $$i \sim j \qua\iff\qua \ga_i(t) = \ga_j(t)\;. $$ 
The partition $\Pi_\al(t)$ refines $\Pi_\al(s)$ if $t>s$. Further, $\Pi_\al(t)$ is 
discrete 
for all sufficiently large $t$; and it is indiscrete when $t$ is sufficiently close to 
zero.  
Hence 
there is a unique affine map $\rho\co[0,1] \to [0,\infty)$ such that  $\Pi_\al(\rho(t))$ 
is 
non-trivial 
exactly for $0<t<1$. 

We can now define the point $\hpi(\al)$ in the nerve $\La_n$ of the partition 
lattice, 
as follows. 
The barycentric coordinate of $\hpi(\al)$ with respect to any non-trivial partition 
$\Pi$ 
is to be the 
length  
of the interval \hbox{$\{t\in [0,1] \;|\; \Pi_\al(\rho(t)) =\Pi\}$}. 
The map $\hpi\co T_n \to \La_n$ so defined is a homeomorphism because the 
fully-grown tree 
$\al$ is uniquely 
determined by its one-parameter family of partitions $\Pi_\al$.  The map is clearly 
equivariant. 
\endprf

\proc{Corollary}\label{2.8} The right action of $\Sgm_n$ on the nerve of the 
partition lattice $\La_n$ extends to a right action of $\Sgm_{n+1}$. 
\end{Corollary}

\proof  This is true of the action on $T_n$, since 
the action of $\Sgm_n$ there is the 
restriction of the action of the group $\Sgm_{n+1}$ of permutations of the full set 
$\{0,1,\dots,n\}$ of labels. \endprf

\section{The representation ${\goths}_n$ and the fundamental 
cycle}\label{3}

\subsection{The Lie representations}

Let ${\cal L}_{n}$ be the free Lie ring on the set of generators $\{x_i\}_{1\le i\le 
n}$. 
We denote by $\Lie_n$ the {\it $n$-linear part \/} of ${\cal L}_n$. This can be 
described 
in many different ways. First, it is  the direct summand of ${\cal L}_n$ spanned by 
all Lie monomials containing each of the $n$ generators exactly once.  
Second, it is isomorphic to the module of all 
natural transformations $\Phi^{\otimes n} \to \Phi$, where $\Phi$ is the forgetful 
functor from Lie 
rings to abelian groups.   Third, it is the $n$th module in the Lie operad.

It is a standard exercise in the use of the Jacobi identity to show that 
$\Lie_n$ is a free abelian group of rank $(n-1)!$, and that the left-regulated Lie 
brackets 
$$\lambda_{\sgm} \qua=\qua [x_{\sgm(1)},[x_{\sgm(2)},[...[x_{\sgm(n-1)},x_n]..]]] 
$$ 
for $\sgm\in\Sgm_{n-1}$
form a $\bbbz$-basis for $\Lie_n$ (see for example \cite{Reu} or \cite{Whi}).  
The~symmetric group $\Sgm_n$ acts on $\Lie_n$ on the left by permuting the $n$ 
generators. 
The~(integral) representation thus afforded is known as the {\it Lie 
representation.\/} 
The subgroup $\Sgm_{n-1}$ 
permutes the above standard basis simply and transitively. Thus the restricted 
representation $\Res^{\Sgm_n}_{\Sgm_{n-1}}\!\Lie_n$ is the regular integral 
representation 
of $\Sgm_{n-1}$.

We denote by $\Lie^*_n$ the dual representation $\Hom(\Lie_n,\bbbz)$. This has 
the 
same character 
as $\Lie_n$, and therefore $\Lie_n$ and $\Lie^*_n$ are isomorphic over the 
rationals. 
Over the 
integers, the two are distinct. 

\proc{Example} $\phantom{99}$\label{3.2}

Take $n=3$, and consider the essentially unique
$\Sgm_3$-invariant Euclidean metric on the two-dimensional space
$\Lie_3 \otimes \,\bbbr$. Both $\Lie_3$ and $\Lie^*_3$ are hexagonal
lattices in this (self-dual) Euclidean space. However, in $\Lie_3$ the
shortest vectors have the form $[x_i,[x_j,x_k]]$, and so these are
$(-1)$-eigenvectors of involutions $(j,k)\in\Sgm_3$.  In $\Lie_3^*$
the shortest vectors are $(+1)$-eigenvectors of involutions.
\end{Example}

Let $\eps$ denote the sign character of $\Sgm_n$. Then we have the
twisted module $\eps\Lie_n$ which is the same abelian group $\Lie_n$
with the twisted action of $\Sgm_n$ defined by
$$\sgm \cdot f(x_1,\ldots,x_n) =
\eps(\sgm)\,f(x_{\sgm(1)},\ldots,x_{\sgm(n)})$$ for every multilinear
Lie monomial $f$ and every $\sgm \in \Sgm_n$.  We also denote the
module $\eps\Lie_n$ by $\eps \otimes \Lie_n$ when we need to stress
that it is the tensor product with the sign module.  We similarly have
the twisted representation $\eps\Lie_n^*$, defined analogously.

\subsection{Connection with Lie superrings}

There is a reinterpretation of the module $\eps\Lie_n$ which significantly 
simplifies 
our  theory.  In place of the free Lie ring on the set of generators 
$\{x_i\}_{1\le i\le n}$ we 
may consider the free Lie superring (or Lie superalgebra over $\bbbz$) 
having the $x_i$ as odd-degree generators.  This we denote by~${\cal S}{}_n$.  
It is a graded algebra 
${\cal S}{}_n^{\;\! even} \oplus {\cal S}{}_n^{\;\! odd}$ 
freely generated by the $x_i$ subject to the relations 
$$\begin{aligned}\relax
[a, b] \qua &=\qua (-1)^{|a||b|+1}[b, a] \\
[[a, b], c] \qua &= \qua [a, [b, c]] - (-1)^{|a||b|}[b, [a, c]] 
\end{aligned} $$
when $a$ and $b$ are homogeneous elements of degrees $|a|$ and $|b|$ 
respectively. 

Let ${\goths}{}_n$ denote the {\it super-Lie representation\/} of $\Sgm_n$.  
This we define to be the \hbox{$n$-linear} part of ${\cal S}_n$: it is spanned by 
the super-Lie monomials containing each $x_i$ exactly once, and 
is a subgroup of ${\cal S}{}_n^{\;\! even}$ or ${\cal S}{}_n^{\;\! odd}$ 
according to the parity of~$n$.  The symmetric group acts by permuting the 
generators, exactly as in the Lie representation.  
We now determine the structure of ${\goths}_n$. 

\proc{Proposition}\label{3.4} The super-Lie representation ${\goths}{}_n$ and the 
twisted 
Lie representation $\eps\Lie_n$ are isomorphic $\,\Sgm_n$-modules. 
\end{Proposition}

\proof  The module $\Lie_n$ is spanned by all Lie 
monomials $g(x_1,x_2,\dots,x_n)$ in 
which each of the $n$ generators occurs once.  
The $n$ variables occur in $g$ in a 
definite order: let this order, reading from left to right, 
be $x_{\ga(1)}, x_{\ga(2)},\ldots, x_{\ga(n)}$.  This gives a permutation 
$\ga \in \Sgm_n$ depending upon~$g$.  
Every such $g$ is the bracket $[h,k]$ of two monomials 
in complementary subsets of the variables.

Let $\calj$ be any Lie superring, and $y_1,y_2,\ldots,y_n$ any elements 
of odd degree in $\calj$.  We claim that there is a unique homomorphism 
$\Tha\,\co \Lie_n \to \calj$ of abelian groups such that for each monomial 
$g(x_1,x_2,\dots,x_n)$ as above 
$$\Tha(g(x_1,x_2,\dots,x_n)) \qua = 
\qua \eps(\ga) G(y_1,y_2,\dots,y_n)$$ where 
$G(y_1,y_2,\dots,y_n)$ is the element of $\calj$ obtained by 
replacing all the Lie bracket operators in $g$ by Lie superbrackets, then 
substituting $y_i$ for $x_i$; and $\eps(\ga)$ is the sign of $\ga(g)$.  The 
uniqueness 
of $\Tha$ is clear. 

We construct $\Tha$ by induction on $n$.  Suppose  the result is known 
already for fewer than $n$ generators.  Then if $g$ has $n$ variables 
and $g = [h,k]$ as above, we have 
$\Tha(g) = \eta\cdot[\Tha(h),\Tha(k)]$ where $[\Tha(h),\Tha(k)]$ 
is the superbracket and $\eta$~is the sign of the inverse shuffle which 
places the generators occurring in $h$ before those occurring in $k$. 
In order to show that $\Tha$ is a well-defined homomorphism in the 
$n\;\!$-variables case, we must verify that it respects the antisymmetry and 
Jacobi relations $[h,k] = -\,[k,h]$ and $[h,[k,l]] + [k,[l,h]] + [l,[h,k]] = 0$ 
of $\Lie_n$.  But this is true because the shuffle signs 
convert these relators into the corresponding signed relators, which 
are zero in the Lie superring $\calj$. The induction is now complete.  

It is not at all difficult to give an instance where $\Tha$ is a 
monomorphism: one begins by taking the $y_i$ to be elementary 
matrices representing elements of the endomorphism superring of a 
superspace of large rank.  
Therefore $\Tha$ is also a monomorphism when $\calj$ is the free 
superring generated by the odd-degree elements $y_1,y_2,\ldots,y_n$. 
So in this case $\Tha$ is an isomorphism of abelian groups onto 
its image, which is ${\goths}{}_n$.  
For any $\sgm\in \Sgm_n$ we have $$\Tha(\sgm\cdot g) = 
\eps(\sgm\ga) (\sgm\cdot G) = \eps(\sgm)(\sgm\cdot\Tha(g))\;.$$
Hence $\Tha$ induces an 
isomorphism between $\Lie_n$ and $\eps\otimes{\goths}{}_n$, and 
between $\eps \Lie_n$ and~${\goths}{}_n$. \endprf

The construction of  the isomorphism between 
$\eps \Lie_n$ and ${\goths}{}_n$ makes use of the 
natural ordering of the set $x_1,x_2,\dots,x_n$ of generators.  The effect of 
adopting a different fundamental ordering, related to the natural one by a 
permutation $\tau$, is to multiply the isomorphism by the sign of $\tau$. 
When the set of generators is unordered, the isomorphism 
is uniquely defined only up to sign. 

\proc{Corollary}\label{3.5} The left-regulated Lie superbrackets 
$$\lambda_{\sgm} \qua=\qua 
[x_{\sgm(1)},[x_{\sgm(2)},[...[x_{\sgm(n-1)},x_n]..]]]\;, $$
where $\sgm\in\Sgm_{n-1}$, form a basis for ${\goths}{}_n$. 
\end{Corollary}
\proof  This follows from \ref{3.1} and \ref{3.4}, or by direct computation as for 
$\Lie_n$. \endprf

The first advantage of ${\goths}{}_n$ over $\eps\Lie_n$ in our homology 
theory is that the Lie superbracket is extremely useful in building the 
fundamental orientation cycle.  
The second advantage is that the isomorphism of ${\goths}{}_n$ with the 
cohomology of $T_n$, which we shall construct in Theorem \ref{4.1} below, 
is totally natural: it does not depend upon an ordering of the set of 
non-zero labels.

\subsection{Construction of the fundamental cycle}

We now describe the construction of a certain cycle $F_n$ in the top-dimensional 
simplicial chain group $C_{n-3}(T_n;{\goths}{}_n)$, which after \ref{3.4} 
we know to be  isomorphic 
to $C_{n-3}(T_n;\eps\Lie_n)$. A top-dimensional simplex $X$ of 
$T_n$ is an $n$-tree which is {\it binary\/}; that is, every internal node is the 
meet of three edges (see \S\ref{2}).  The vertices 
of the simplex correspond to the internal edges of the tree $X$; 
so an orientation of the 
simplex may be specified by a word $w_X = e_0e_1\ldots e_{n-3}$ in which 
each internal edge of $X$ appears once. We write this oriented simplex~
$\lan w_X \ran$. If~the order of the edges in $w_X$ is altered by a 
permutation $\sgm$, then $\lan w_X \ran$ is multiplied by $\eps(\sgm)$.  
We shall choose such an ordering for each $(n-3)$-simplex $X$ of $T_n$, and 
we shall set 
$$F_n\qua = \qua \sum_X \;\lan w_X \ran\otimes c_X \qua\in\qua 
C_{n-3}(T_n;{\goths}_n)$$ 
where $c_X\in {\goths}_n$ is a coefficient monomial which we now describe. 

We use induction on $n$ to define $w_X$ and $c_X$. When $n\le 2$ the complex 
$T_n$ is empty, as the 
unique $n$-trees $X$ have no internal edges; but nevertheless we formally define 
$c_X = x_1$ when $n=1$ and $c_X = - [x_1,x_2]$ when $n=2$.  Here, and 
throughout 
the rest of \S\ref{3}, the bracket is the operation in a free Lie {\it superring \/}.  
The variables 
$x_i$ are in bijective correspondence with the non-zero tree labels. 
Since the $x_i$ are assigned odd grading, we have $[x_1,x_2]= [x_2,x_1]$.  
The total ordering on the set of labels has no significance.  

We continue with our inductive definition. Suppose now that $n>2$.   

The $n$-tree $X$ is obtained by grafting together at the root 
an $i$-tree $Y$ and an $(n-i)$-tree $Z$, where $0<i<n$,\/ 
as in the following diagram 

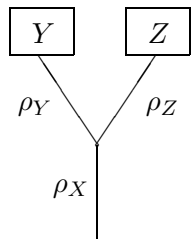
\setlength{\unitlength}{0.44pt}
\begin{picture}(450,250)(-400,-100)

\def\framebox(#1,#2){\put(0,0){\line(0,1){#2}} \put(0,0){\line(1,0){#1}}
\put(0,#2){\line(1,0){#1}} \put(#1,0){\line(0,1){#2}}}

\put(0,0){\line(0,-2){80}}
\put(0.3,2){\line(2,3){50}}
\put(-0.3,2){\line(-2,3){50}}
\put(0,0){\circle{2}}

\put(-70,30){$\rho_Y$}
\put(40,30){$\rho_Z$}
\put(-40,-40){$\rho_X$}
\put(25,78){\framebox(55,40)}
\put(42,88){$Z$}
\put(-75,78){\framebox(55,40)}
\put(-58,88){$Y$}
\end{picture}

since the root of $X$ meets two other edges. 
The internal edges of $X$ comprise the internal edges of $Y$ and those 
of~$Z$, together with 
the roots $\rho_Y$ and $\rho_Z$ of $Y$ and $Z$. Suppose by induction that 
we have well-defined terms $\lan w_Y\ran \otimes c_Y$ and 
$\lan w_Z\ran \otimes c_Z$ corresponding to orderings of 
the internal edges of $Y$ and $Z$, the coefficients $c_Y$ and $c_Z$ being 
super-Lie 
monomials in the sets of non-zero labels of $Y$ and $Z$ respectively.  
We set $w_X =  \rho_Yw_Y\rho_Zw_Z$, this being an ordering 
of all the internal edges of $X$.  
The set of non-zero labels of $X$ is the disjoint union of those of 
$Y$ and $Z$. We put $c_X = (-1)^{|Y|}[c_Y,c_Z]$ where 
$|Y|$ is the number $i$ of leaves, excluding the root, in $Y$.  This is a 
super-Lie monomial in the non-zero labels of $X$. 

To justify this construction, we verify in the next two lemmas that 
$\lan w_X\ran \otimes c_X$ is a well-defined element 
of the chain group $C_{n-3}(T_n;{\goths}_n)$, and that the sum $F_n$ of these 
elements over all top-dimensional simplices $X$ is a cycle. 

\proc{Lemma}\label{3.7} The term $$\lan w_X \ran\otimes c_X \qua\in\qua 
C_{n-3}(T_n;{\goths}_n)$$ 
is well defined for every binary $n$-tree $X$, and is $\Sgm_n$-invariant. 
\end{Lemma}

\proof  We again use induction on $n$. For $n=1$ the result is completely 
trivial. 
When 
$n\ge 2$ the tree $X$ is obtained by grafting to the root an $i$-tree $Y$ and an 
$(n-i)$-tree 
$Z$, where $0<i<n$ and therefore both $\lan w_Y \ran\otimes c_Y$ and 
$\lan w_Z \ran\otimes c_Z$ 
are well-defined by inductive hypothesis. 

Since the oriented simplices $\lan w_Y \ran$ and $\lan w_Z \ran$ are evidently 
well defined up to sign, the coefficients $c_Y$ and $c_Z$ are 
well defined up to sign.  Therefore the total sign in the formula 
$$\lan w_X \ran\otimes c_X\qua = \qua \lan \rho_Yw_Y\rho_Zw_Z \ran \otimes 
(-1)^{|Y|}[c_Y,c_Z]$$ 
is independent of the orientations of the two simplices $Y$ and $Z$, by bilinearity 
of tensor product. 
The only remaining ambiguity is that the sub-trees $Y$ and $Z$ can be listed in 
either order. 
Interchanging them may affect the orientation of the simplex $X$, and may alter  
the sign of $c_X = (-1)^{|Y|}[c_Y,c_Z]$. 

We compute the change in the total sign when $Y$ and $Z$ are interchanged. 
Since the number 
of internal edges in a binary tree equals the number of leaves (excluding the root) 
less two, the effect of interchanging 
$\rho_Yw_Y$ and $\rho_Zw_Z$ is to multiply $\lan w_X\ran$ by the sign 
$(-1)^{(i-1)(n-i-1)}$. 
On the other hand, the replacement of $(-1)^{|Y|}$ by $(-1)^{|Z|}$ introduces a 
sign $(-1)^n$, 
whilst the replacement of the superbracket $[c_Y,c_Z]$ by $[c_Z,c_Y]$ gives a 
factor 
$(-1)^{i(n-i)+1}$, because $c_Y$ and $c_Z$ have degrees $i$ and $n-i$ 
respectively.  The product of these four 
signs is $+1$, so that $\lan w_X \ran\otimes c_X$ is independent of the order in 
which $Y$ and $Z$ are taken. The inductive step is thus complete. 
Finally, the $\Sgm_n$-invariance of $\lan w_X \ran \otimes c_X$ is obvious, since 
the ordering of the set of labels has nowhere been used. The lemma is proved. 
\endprf

\proc{Proposition}\label{3.8} The $(n-3)$-chain 
$$F_n\qua = \qua \sum_X \;\lan w_X \ran\otimes c_X \qua\in\qua 
C_{n-3}(T_n;{\goths}_n)$$ 
is a $\Sgm_n$-invariant cycle. \end{Proposition}

\proof  The $\Sgm_n$-invariance of $F_n$ follows from Lemma \ref{3.7}. 

To prove that $F_n$ is a cycle, we calculate all the coefficients of the 
$(n-4)$-chain 
$\dee F_n$. Every oriented $(n-4)$-simplex $\sgm$ corresponds to a tree $S$ 
which is binary 
apart from one exceptional node, where four edges meet. There are two cases, 
depending upon the position of the exceptional node. 

Suppose first that the exceptional node is the root node. Then $S$ is the result of 
grafting 
together three trees $X$, $Y$ and $Z$, as shown on the left here. As orientation 
for $\sgm$ we may take the edge-sequence 
$\lan \rho_X w_X\rho_Yw_Y\rho_Zw_Z\ran$, since this lists all the internal edges 
(in the notation of \ref{3.6}). One of the three $(n-3)$-simplices incident upon $\sgm$ is 
 obtained by grafting $X$ to $Y$ (which introduces a new root $\rho_{XY}$) and 
grafting the result to $Z$, as in the picture on the right. 

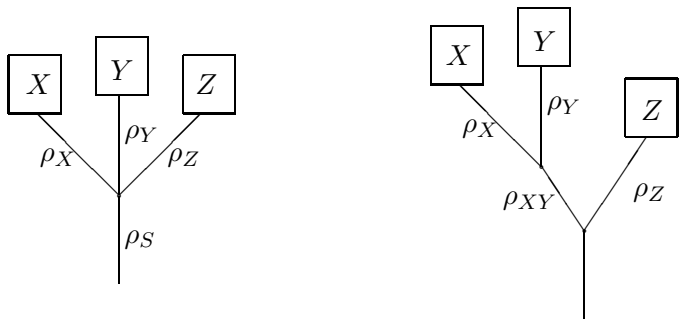
\setlength{\unitlength}{0.44pt}
\begin{picture}(1200,310)(-400,-120)

\put(-600,-100){\setlength{\unitlength}{0.44pt}
\begin{picture}(450,250)(-400,-100)

\def\framebox(#1,#2){\put(0,0){\line(0,1){#2}} \put(0,0){\line(1,0){#1}}
\put(0,#2){\line(1,0){#1}} \put(#1,0){\line(0,1){#2}}}

\put(0,0){\line(0,-1){75}}
\put(0.3,1){\line(1,1){70}}
\put(-0.3,1){\line(-1,1){70}}
\put(0,1){\line(0,1){85}}
\put(0,0){\circle{2}}

\put(-70,30){$\rho_X$}
\put(40,30){$\rho_Z$}
\put(3,-40){$\rho_S$}
\put(3,50){$\rho_Y$}

\put(55,71){\framebox(45,50)}
\put(64,88){$Z$}

\put(-20,87){\framebox(45,50)}
\put(-9,100){$Y$}
\put(-95,71){\framebox(45,50)}
\put(-82,88){$X$}
\end{picture}}

\put(-200,-130){\setlength{\unitlength}{0.44pt}
\begin{picture}(450,250)(-400,-100)

\def\framebox(#1,#2){\put(0,0){\line(0,1){#2}} \put(0,0){\line(1,0){#1}}
\put(0,#2){\line(1,0){#1}} \put(#1,0){\line(0,1){#2}}}

\put(0,0){\line(0,-1){75}}
\put(0.3,1){\line(2,3){53}}
\put(-0.3,1){\line(-2,3){36}}
\put(0,0){\circle{2}}

\def\subtree{
\put(-0.5,1){\line(-1,1){70}} \put(0,1){\line(0,1){85}} 
\put(0,0){\circle{2}}
\put(-70,30){$\rho_X$} \put(3,50){$\rho_Y$}
\put(-20,87){\framebox(45,50)}
\put(-9,100){$Y$}
\put(-95,71){\framebox(45,50)}
\put(-84,88){$X$}
}

\put(-37,55){\subtree}

\put(-72,23){$\rho_{XY}$}

\put(40,30){$\rho_Z$}

\put(35,81){\framebox(45,50)}
\put(47,96){$Z$}

\end{picture}}

\end{picture}

This simplex can be oriented as $\lan 
\rho_{XY}\rho_Xw_X\rho_Yw_Y\rho_Zw_Z\ran$: 
then its incidence
index upon $\sgm$ is $+1$, and its contribution 
to the coefficient of $\dee F_n$ 
on $\sgm$ is, by 
two applications of the construction in \ref{3.6} 
$$(-1)^{|Y|}[[c_X,c_Y],c_Z].$$ 
\noindent The other two $(n-3)$-simplices incident on $\sgm$ 
are obtained from the first by permuting 
$X$, $Y$ and $Z$ cyclically. The simplex $\lan 
\rho_{YZ}\rho_Yw_Y\rho_Zw_Z\rho_Xw_X\ran$ has 
incidence index $(-1)^{(|X|+1)(|Y|+|Z|)}$ upon $\sgm$ (oriented as before), and 
the applications of the construction of \ref{3.6} 
contribute altogether a sign $(-1)^{|Z|}$. 
The second simplex therefore contributes 
$$(-1)^{|Y| + |X||Y| + |X||Z|}[[c_Y,c_Z],c_X]. $$
By cyclic permutation, the third simplex incident upon $\sgm$ contributes 
$$(-1)^{|Y| + |Z||X| + |Z||Y|}[[c_Z,c_X],c_Y]. $$ 
The total coefficient of $\dee F_n$ on $\sgm$ is thus zero 
in view of the Jacobi identity 
$$(-1)^{|Y||Z|}[[c_Z,c_X],c_Y] + (-1)^{|Z||X|}[[c_X,c_Y],c_Z] + 
(-1)^{|X||Y|}[[c_Y,c_Z],c_X] 
= 0$$
in the Lie superring.  The first case is proved. 

In the second case, the exceptional node (where four edges meet) of $S$ is not 
the 
root node. 
Thus $S$ is obtained by grafting together at the root a binary tree $T$ 
and a 
tree $U$ 
having one exceptional node. Let $m$ be the number $|U|$ of leaves of $U$. 
It follows by application of the definition in \ref{3.6} that 
the coefficient of $\dee F_n$ on $\sgm$ is $[c_T, e]$ where $e$ is the coefficient 
of $\dee F_m$ on an appropriate orientation of the simplex represented by $U$. 
By induction on $n$, we know that $F_m$ is a cycle, so $e$ is zero. 
This establishes the second case, and completes the proof of~\ref{3.8}. \endprf 


\proc{Example}\label{3.9} $\phantom{99}$

We illustrate the foregoing theory by calculating the 
fundamental cycle $F_5$  of~$T_5$.  
To reinforce the point that the ordering of the tree labels 
is nowhere used, we label the trees by $\{0,a,b,c,d,e\}$ in place of $\{0,1,2,3,4,5\}$. 
We use $a,b,c,d,e$ also to denote the corresponding odd-degree generators of 
the Lie superring. 

We enumerate the top-dimensional simplices of $T_5$.  There are 105 of these, 
generated by the trees pictured here and the action of the symmetric group 
$\Sgm_5$ on the labels $a,b,c,d,e$.  We choose orientations for the simplices 
by ordering the internal edges, and then calculate the 105 corresponding terms 
in the cycle~$F_5$.  Each term is obtained by using the formula of \ref{3.6} twice.

There are 60 terms of the form 
$\Phi = \lan x,y,z \ran \otimes [a, [b, [c, [d, e]]]] $ 
corresponding to trees of the shape:

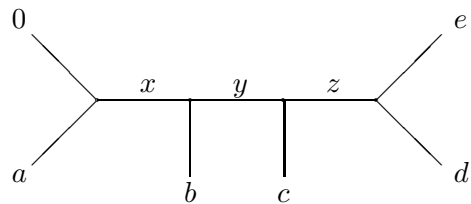
\setlength{\unitlength}{0.44pt}
\begin{picture}(450,230)(-380,-120)

\put(-121,-1){\line(-1,-1){55}}
\put(-121,1){\line(-1,1){55}}
\put(-120,0){\line(1,0){80}} 
\put(-40,0){\line(1,0){80}}
\put(40,0){\line(1,0){80}}
\put(121,1){\line(1,1){55}}
\put(121,-1){\line(1,-1){55}}

\put(-40,-0.5){\line(0,-1){65}}
\put(41,-0.5){\line(0,-1){65}}

\put(-120,0){\circle{2}}
\put(-40,0){\circle{2}}
\put(40,0){\circle{2}}
\put(120,0){\circle{2}}


\put(-85,7){$x$}
\put(-5,7){$y$}
\put(75,7){$z$}
\put(-195,62){$0$}
\put(-195,-70){$a$}
\put(-47,-88){$b$}
\put(33,-88){$c$}
\put(185,62){$e$}
\put(185,-70){$d$}
\end{picture}

another 30 terms of the form $\Psi = - \lan p,q,r \ran \otimes [[a,
b], [c, [d, e]]]$ corresponding to trees of the shape:

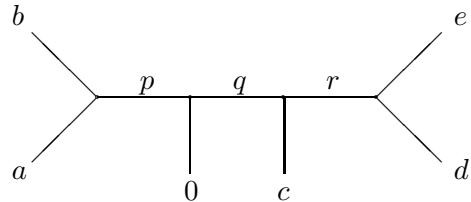
\setlength{\unitlength}{0.44pt}
\begin{picture}(450,230)(-380,-120)

\put(-121,-1){\line(-1,-1){55}}
\put(-121,1){\line(-1,1){55}}
\put(-120,0){\line(1,0){80}} 
\put(-40,0){\line(1,0){80}}
\put(40,0){\line(1,0){80}}
\put(121,1){\line(1,1){55}}
\put(121,-1){\line(1,-1){55}}

\put(-40,-0.5){\line(0,-1){65}}
\put(41,-0.5){\line(0,-1){65}}

\put(-120,0){\circle{2}}
\put(-40,0){\circle{2}}
\put(40,0){\circle{2}}
\put(120,0){\circle{2}}


\put(-85,7){$p$}
\put(-5,7){$q$}
\put(75,7){$r$}
\put(-195,62){$b$}
\put(-195,-70){$a$}
\put(-47,-88){$0$}
\put(33,-88){$c$}
\put(185,62){$e$}
\put(185,-70){$d$}

\end{picture}

and 15 terms of the form 
$\Om = - \lan i,j,k \ran \otimes [a, [[b, c], [d, e]]]$ corresponding to trees of the shape:

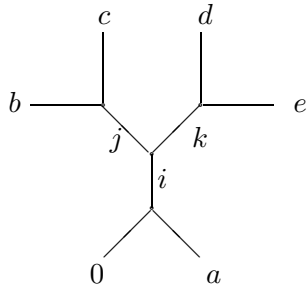
\setlength{\unitlength}{0.44pt}
\begin{picture}(450,280)(-380,-140)

\put(-3,1){\line(-1,1){40}}
\put(1,1){\line(1,1){40}}
\put(0,-1){\line(0,-1){45}} 

\put(44,42){\line(1,0){60}}
\put(42,44){\line(0,1){60}}
\put(-44,42){\line(-1,0){60}}
\put(-42,44){\line(0,1){60}}
\put(-1,-49){\line(-1,-1){40}}
\put(1,-49){\line(1,-1){40}}

\put(0,0){\circle{2}}
\put(-42,42){\circle{2}}
\put(42,42){\circle{2}}
\put(0,-47){\circle{2}}


\put(-38,6){$j$}
\put(33,6){$k$}
\put(3,-30){$i$}

\put(-125,36){$b$}
\put(-48,111){$c$}
\put(-55,-111){$0$}
\put(45,-111){$a$}
\put(120,36){$e$}
\put(38,111){$d$}


\end{picture}

The number of terms in each case is given by dividing $5\>\!!= 120$ by 
the order of the symmetry group of the tree fixing the root. 
Our construction takes care of the signs.   
We have $$\begin{aligned}F_5 \qua = 
\qua \sum \;\lan x,y,z \ran \otimes [a, [b,\ &[c, [d, e]]]] 
		\qua - \qua \sum \; \lan p,q,r \ran \otimes [[a, b], [c, [d, e]]]  \\
- \qua &\sum \;\lan i,j,k \ran \otimes [a, [[b, c], [d, e]]]
\end{aligned} $$
where there are 105 terms in all, corresponding to all the labelled
trees up to label-preserving isomorphism.  (The letters attached to the
internal edges are bound variables over which summation is carried
out: other symbols can be substituted for them at will.)  Now we check
that the boundary $\dee F_5$ is zero.  We revert to the notation
$\Phi, \Psi, \Om$ for the three types of terms in $F_5$.  On each type
there are three face operators to consider.

First we look at the terms $\dee_2 \Phi$.  We have  $\dee_2 \Phi = \lan x,y \ran 
\otimes [[a, b], [c, [d, e]]]$ where $\lan x,y \ran$ 
is an orientation of the tree:

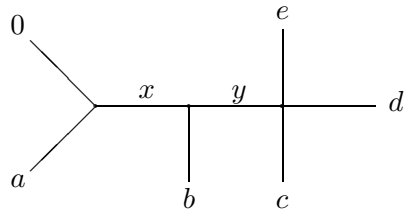
\setlength{\unitlength}{0.44pt}
\begin{picture}(450,250)(-450,-120)

\put(-121,-1){\line(-1,-1){55}}
\put(-121,1){\line(-1,1){55}}
\put(-120,0){\line(1,0){80}} 
\put(-40,0){\line(1,0){80}}
\put(40,0){\line(1,0){80}}

\put(-40,-0.5){\line(0,-1){65}}
\put(41,-0.5){\line(0,-1){65}}
\put(41,0.5){\line(0,1){65}}

\put(-120,0){\circle{2}}
\put(-40,0){\circle{2}}
\put(40,0){\circle{2}}


\put(-85,7){$x$}
\put(-5,7){$y$}

\put(-195,62){$0$}
\put(-195,-70){$a$}
\put(-47,-88){$b$}
\put(33,-88){$c$}
\put(33,75){$e$}
\put(130,-4){$d$}
\end{picture}

in which the leaves $c, d$ and $e$ are attached at the same node.  Therefore the 
60 terms $\dee_2 \Phi$ cancel among themselves in threes by virtue of the 
relation 
$$[c, [d, e]] \; + \; [d, [e, c]] \; + \; [e, [c, d]] \qua = \qua 0 $$ 
in the coefficient module ${\goths}_5$.  In the same way, the 30 terms 
$\dee_2\Psi$ 
have sum zero, and the 60 terms $\dee_0 \Phi$ cancel the 30 terms $\dee_0\Psi$ 
by the supercommutator formula 
$$[[a, b], [c, [d, e]]] \; - \; [a, [b, [c, [d, e]]]] \; - \; [b, [a, [c, [d, e]]]]\qua = \qua 0 \;.$$ 
Finally each of the 15 terms $\dee_0\Om$ cancels with a pair of terms 
$\dee_1\Psi$ through 
$$[a, [[b, c], [d, e]]]\qua = \qua [[b, c], [a, [d, e]]]\; - \; [[d, e], [a, [b, c]]] $$ 
and each of the terms in $\dee_1\Om$ or $\dee_2\Om$ likewise cancels a pair 
of terms in $\dee_1\Phi$.  This accounts for all 315 terms in $\dee F_5$, and 
verifies that $F_5$ is a cycle. \end{Example}

\section{The homology and cohomology as integral 
representations of  $\Sgm_n $}\label{4}

We showed in \S\ref{2} that the space $T_n$ has reduced homology only in dimension 
$n-3$. 
We now determine the action of $\Sgm_n$ on this finite $\bbbz$-module.  The 
following theorem 
is equivalent to a result proved by M. Wachs \cite{Wac} for the cohomology of the 
poset of partitions. 

\proc{Theorem}\label{4.1} There are isomorphisms of left $\,\Sgm_n$-modules 
and of right $\,\Sgm_n$-modules respectively 
$$\begin{aligned}
\til H^*(T_n;\bbbz)\qua &\approx \qua {\goths}_n \qua \approx
\qua\eps \,{\Lie}_n \\
\til H_*(T_n;\bbbz)\qua &\approx \qua {\goths}_n^{*} \qua 
\approx \qua \eps\,{\Lie}_n^*.
\end{aligned} $$
\end{Theorem}

\proof  As $T_n$ is a finite complex with torsion-free homology, there is for 
every 
coefficient group $G$ a natural universal coefficient isomorphism  
$$ \til H_*(T_n;\,G) \qua\approx\qua {\Hom}_{\bbbz}(\til 
H^*(T_n),\,G)\;.$$
In particular, the integral homology and integral cohomology  are dual 
$\bbbz$-mod\-ules. 
In the above we take $G$ to be ${\goths}_n$ or the isomorphic module 
$\eps\Lie_n$. 
The homology class of the invariant cycle $F_n$ constructed 
in \ref{3.6} corresponds under the universal coefficient 
isomorphism to some homomorphism 
\hbox{$\theta_n\co \til H^*(T_n) \to {\goths}_n$}, which is indeed  
$\Sgm_n$-equivariant because 
the universal coefficient isomorphism is natural and $F_n$ is $\Sgm_n$-invariant. 

Let $f_{\sgm}$ be the cocycle of $T_n$ which takes the value $1$ on the simplex 
$\ga_{\sgm}$ depicted in \ref{2.5}, and is zero on all other $(n-3)$-simplices of $T_n$. 
Then the set $$\{[f_{\sgm}]\mid \sgm \in \Sgm_{n-1}\}$$ is a basis 
for $\til H^*(T_n)$, by \ref{2.5}. 
The value of $\tha_n$ on $[f_{\sgm}]$ is $\pm \la_{\sgm}$, 
by \ref{3.6} and induction on $n$, where $\la_{\sgm}$ is the basis element 
$[x_{\sgm(1)},[x_{\sgm(2)},[...[x_{\sgm(n-1)},x_n]..]]]$ of the 
$\bbbz$-module ${\goths}_n$, as in \ref{3.5}. 
Therefore we have shown 
that the $\Sgm_n$-homomorphism $\tha_n$ carries a basis for $\til H^*(T_n)$ 
into a basis for ${\goths}_n$, which is isomorphic to $\eps\Lie_n$ by \ref{3.4}. 
That proves the first assertion in the theorem. The second assertion 
follows by taking $\bbbz$-duals. \endprf 

\proc{Corollary}\label{4.2} The $\Sgm_n$-module $\Lie_n$ is the restriction of a left 
$\,\Sgm_{n+1}$-module. \end{Corollary}

\proof  We have seen that $\Sgm_{n+1}$ acts on $T_n$, so it acts upon 
$\eps\til H^*(T_n)$, which is isomorphic to $\Lie_n$. \endprf 

\proc{Definition}\label{4.3} \rm We denote by $\hatLie_n$ the left $\Sgm_{n+1}$-module 
so defined, 
and by  $\hatLie_n^*$ the dual right $\Sgm_{n+1}$-module. 
Evidently $\Res^{\Sgm_{n+1}\mathstrut}_{\Sgm_n}\hatLie_n\;\approx\; {\Lie\>\!}_n
$.  
\end{Definition}

\section{The Whitehouse extension}\label{5}

Let $H$ be a subgroup of the finite group $G$. 
Since the functor $\Ind_H^G$ is left adjoint to $\Res_H^G$, there is an adjunction 
morphism 
$\Ind_H^G \Res_H^G M \lra M$ for every $G$-module $M$, and it is always 
surjective. 

The following theorem was proved algebraically by Whitehouse \cite{Whi};
it also follows from work of Sundaram \cite{Sun}. 

\proc{Theorem}\label{5.1} There is a short exact sequence of 
$\Sgm_{n+1}$-modules 
$$ 0 \lra {\Lie}^*_{n+1} \lra  {\Ind}^{\Sgm_{n+1}}_{\Sgm_n}{\Lie}^*_n \lra 
\hatLie_n^* \lra 0 $$ 
in which the surjection is the adjunction morphism connecting 
$\,\Ind^{\Sgm_{n+1}}_{\Sgm_n}\!$ and $\,\Res^{\Sgm_{n+1}}_{\Sgm_n}\!$. 
\end{Theorem}

We shall give a topological proof.  We recall that there is an analogue for 
$G$-spaces of the 
induction functor $\Ind_H^G$, namely smash product with $(G/H)\!_+$.  
Theorem~\ref{5.1} follows immediately from the cofibration in the 
following theorem, by taking reduced homology. 

\proc{Theorem}\label{5.2} There is a $\Sgm_{n+1}$-equivariant cofibration 
$$X \qua \subset \qua T_{n+1} \qua \lra \qua  
(\Sgm_{n+1}/\Sgm_n)_+\wedge 
ST_n $$ 
where $X$ is $\Sgm_{n+1}$-homotopy equivalent to $T_n$ and $S$ denotes 
suspension. \end{Theorem}

\proof  In the tree space $T_{n+1}$ we consider for each $1\le i \le n+1$ 
the star $\st (v_{0i})$ of the vertex $v_{0i}$ of \ref{2.2}.  This consists of all the 
$(n+1)$-trees in which the leaves labelled $0$ and $i$ 
are separated from the rest of the tree by an internal edge of positive length.  
These $n+1$
open stars are therefore disjoint.  The closure $\ost (v_{0i})$ is isomorphic to the 
contractible 
space $\tilde T_n$, by the homeomorphism which deletes the leaf $0$ and then 
relabels 
the leaf $n+1$ as $0$.  Let $X$ be the complement 
$T_{n+1} \;\setminus \;\bigcup_{i=1}^{n+1}\st (v_{0i})$. 
By~\ref{2.3} the subcomplex $\tilde T_n$ is the cone on $T_n$, so 
$$T_{n+1}/X \qua \approx \qua \bigvee_{i=1}^{n+1}\ost (v_{0i})/\dee(\ost 
(v_{0i}))
\qua \approx \qua \bigvee_{i=1}^{n+1} ST_n \;.$$
The isomorphism $\ost (v_{0,n+1})/\dee(\ost (v_{0,n+1}))\approx ST_n$ is 
$\Sgm_n$-equivariant, 
and $\Sgm_{n+1}$ induces a bijection between the summands and the left cosets 
$\Sgm_{n+1}/\Sgm_n$, 
so the wedge on the right above is equivariantly isomorphic to 
$(\Sgm_{n+1}/\Sgm_n)_+\wedge ST_n $. 

We must now prove that $X$ is $\Sgm_{n+1}$-homotopy-equivalent to $T_n$, 
where 
$\Sgm_{n+1}$ 
acts on $T_n$ by permuting the labels $\{0,1,\dots,n\}$, with $0$ in place of 
$n+1$.  
We assume $n \ge 3$: the reasoning below requires trivial modifications for 
$n=1,2$. 

Given any tree in $T_n$, we obtain a tree in $X \subset T_{n+1}$ by the following 
procedure: relabel 
the leaf $0$ as ${n+1}$, then attach a new leaf labelled $0$ at any point on any 
{\it 
internal \/} edge. 
Every point of $X$ arises once only in this way, so the reverse procedure 
(deleting the leaf $0$ and 
relabelling the leaf $n+1$ as $0$) gives a projection $\hpi\co X\to T_n$.  
As usual, when the leaf-deletion results in the amalgamation 
of two edges with lengths $x$ and $y$, the new combined edge is assigned 
length $\max(x,y)$.  

The inverse image under $\hpi$ of any point $t\in T_n$ is isomorphic to the 
subtree consisting of the 
internal edges of $t$.  This is contractible.  By the Vietoris Theorem~\ref{5.3} 
below, $\hpi$ is a homotopy equivalence.  

The map $\hpi$ is $\Sgm_{n+1}$-equivariant by construction.  To show that it is a 
equivariant homotopy 
equivalence, it suffices to show that, for every subgroup $H$ of $\Sgm_{n+1}$, the 
restriction 
$\hpi^H\co X^H\to T^H_n$ between the fixed point sets is a homotopy 
equivalence.  
If 
$t \in T^H_n$ is any $H$-fixed 
point, then the group $H$ acts by isometries on the tree~$t$, permuting the labels. 
 
Whenever a finite 
group acts by isometries on a finite tree, the set of invariant points is contractible.  
In our case the invariant 
points form the inverse image of the point $t$ under $\hpi^H$.  Thus $\hpi^H$ 
satisfies the hypotheses of 
Theorem~\ref{5.3}, and the induced map of fixed point sets is always a homotopy 
equivalence, which proves Theorem~\ref{5.2} and therefore Theorem~\ref{5.1}.  
\endprf 

The Vietoris theorem used in the above proof was proved in the following 
strong form by M. M. Cohen  (\cite{CoM}, Theorem 11.1). 
\proc{Theorem}\label{5.3} 
 A piecewise linear map between finite polyhedra is a simple homotopy 
equivalence 
if the inverse image of every point is contractible. \qed 
\end{Theorem}

\noindent Homotopy equivalence under these hypotheses is a theorem of Smale 
\cite{Sma}.  The simply connected case already follows from Vietoris \cite{Vie}  
by an application of Whitehead's Theorem. 

\sh{Acknowledgements} The author acknowledges the benefit of  
conversations with Teimuraz Pirashvili, Colin Rourke and Brian Sanderson, and 
valuable commentary from the referee.

\Addresses\recd
\end{document}